\def\ZZ         {{\mathbb Z}}
\def\RR         {{\mathbb R}}
\def\CC         {{\mathbb C}}

\def\QQ         {{\mathbb Q}}
\def\PP         {{\mathbb P}}
\def\NN         {{\mathbb N}}
\def\ZZ         {{\mathbb Z}}

\def\CO         {{\mathcal O}}

\def\cal        {\mathcal}

\documentclass{amsart}
\usepackage{amssymb}
\usepackage{verbatim}
\usepackage{eucal}
\usepackage{enumerate}
\usepackage{tikz}
\usetikzlibrary{decorations.markings,decorations.pathmorphing,cd}
\usetikzlibrary{arrows}
\usetikzlibrary{calc}
\usepackage{pgfplots}
\pgfplotsset{every axis/.append style={
                    axis x line=middle,    
                    axis y line=middle,    
                    axis line style={-,color=blue}, 
                    xlabel={$x$},          
                    ylabel={$y$},          
            }}

\usepackage{amsthm}
\usepackage{amssymb}
\usepackage{enumerate}
\usepackage{graphics}
\usepackage{hyperref}
\newtheorem{thm}{Theorem}[section]
\newtheorem{theorem}{Theorem}[section]
\newtheorem{lemma}[theorem]{Lemma}
\newtheorem{prop}[theorem]{Proposition}

\theoremstyle{definition}
\newtheorem{rem}[theorem]{Remark}
\newtheorem{dfn}[theorem]{Definition}

\newtheorem{example}[theorem]{Example}

\theoremstyle{remark}

\title{Families of singular algebraic varieties that are rationally elliptic spaces.} 
\author{A.Libgober}
\address{Department of Mathematics\\
University of Illinois\\
Chicago, IL 60607}
\email{libgober@uic.edu}
\thanks{DAS and COI statements: No data was collected for this research. Results used in this work are all coming
from the references given in the paper.
Author has no conflict of interests presenting this work
}

\begin{document}
\begin{abstract} We discuss families of hypersurfaces with isolated singularities in projective space with 
the property that 
the sum of the ranks of the rational homotopy and the homology groups is finite. They represent infinitely many distinct homotopy types 
with all hypersurfaces having a nef canonical or anti-canonical class. 
In the appendix we show that such an infinite family of smooth rationally elliptic 3-folds does not exist. 
\end{abstract}
\maketitle

\centerline{\dedicatory{\it  In memory of J.W.Wood}}
\section{Introduction}
A rationally elliptic space (cf. \cite{halperin}, Chap. 32.) is a topological space $X$ that is simply connected 
and satisfies  the finiteness condition:
$${\rm dim}  H_*(X,\QQ) <\infty \ \ \ and \ \  {\rm dim}  \pi_*(X) \otimes \QQ < \infty. $$
This homotopy-theoretic property appears to have geometric significance. 
A classical example of how geometry may affect rationally ellipticity is given by the Bott conjecture, predicting that a smooth manifold admitting a metric of nonnegative sectional curvature is rationally elliptic (cf. \cite{grovehalperin}, (1.5)). Another example is a simple observation that homogeneous spaces 
 are rationally elliptic, giving a large pool of samples (cf. \cite{halperin} Chap. 32). 
 
 Many rationally elliptic manifolds are algebraic (e.g. quotients of reductive groups by a parabolic subgroup). 
 The links of simplest isolated singularities give other examples of rationally elliptic spaces in the algebro-geometric context.  Below 
  we provide explicit equations for {\it singular projective algebraic varieties} with Gorenstein singularities and 
  nef canonical or anticanonical bundles  that are rationally elliptic (and for which Bott's geometric relation is not applicable). 
 We also consider the following question: under what geometric conditions one may expect finiteness of homotopy (or deformation) types of rationally elliptic
 algebraic varieties.  

The main result of this paper is a series of examples (cf. (\ref{kollar}), (\ref{myprocspace}),(\ref{newspace})  that give infinite families of {\it singular} rationally elliptic 
hypersurfaces with isolated singularities for any dimension greater than one, having distinct homotopy types but falling into 
two homotopy types over $\RR$. In the Appendix below (using the results in \cite{amorosbiswas}) 
we show that there are no {\it  smooth rationally elliptic} 3-folds 
with the canonical class being nef (there are finitely many deformation types of Fano 3-folds and the work \cite{amorosbiswas} enumerates 
those that are rationally elliptic). However, there are infinitely many homotopy types of rationally elliptic K\"ahler manifold 
for which neither canonical or anti-canonical class are nef (cf. Theorem \ref{lowdim}(b)). The singular examples in this note have nef canonical or anti-canonical class. An interesting question
is if there are only finitely many homotopy types that are represented by smooth rationally elliptic algebraic varieties with nef canonical or anticanical class, as is the case for smooth varieties in dimensions 2 and 3. Note  that there is an {\it infinite} family of {\it real homotopy types of rationally elliptic {\it smooth} $C^{\infty}$ manifolds} in real dimension 6 (cf. \cite{herrmann}, Prop. 3.8)

In addition to the classification of the low dimensional smooth rationally elliptic algebraic varieties studied in 
\cite{amorosbiswas} (surfaces and threefolds), in \cite{yangsu} the Hodge diamonds of K\"ahler rationally elliptic smooth 4-folds are described.
The latter, based on the work \cite{neisendorfer}, also shows that the only rationally elliptic smooth complete intersections 
are projective spaces and quadrics.

A related problem of classification of K\"ahler manifolds with the simplest homotopy types (and are rationally elliptic) was considered 
in higher dimensions much earlier.
The case of $\ZZ$-homotopy type of $\PP^n$ was studied in 
\cite{wood} where it was shown that for $n \le 6$, $\PP^n$ is the only K\"ahler manifold with its homotopy type (for $n=7$, cf. \cite{debarre} for a partial extension 
and also \cite{liping}). Fixing only $\QQ$-homotopy type allows for additional solutions (for example, some Fano 3-folds; cf. Remark \ref{fano} below). 
The problem addressed in \cite{wood} is a refinement of the problem of classifying K\"ahler manifolds {\it homeomorphic} to 
$\PP^n$ solved by Hirzebruch and Kodaira (cf. \cite{kodaira}) and S.T.Yau (cf.\cite{yau}).  Brieskorn obtained results similar to \cite{kodaira} in the case 
 of the homeomorphism type of a quadric having dimension $\ne 2$. Our examples (\ref{newspace}) explicitly give {\it the singular spaces} having the same homology as 
 a smooth quadric, although homotopically equivalent to quadrics only over $\RR$ (cf. Prop. \ref{typeofquadric} for a precise statement). 



Homological projective spaces and quadrics arise naturally in the study of the Montgomery-Yang problem on 
circle actions (cf. \cite{kollar5manifolds} for a review). Singular hypersurfaces that are homological projective spaces 
lead to circle actions on stratified spaces that are homology spheres, and  we show that J.Koll\'ar's examples in \cite{kollaryau} also 
serve as examples of rationally elliptic spaces. Similarly, our examples (\ref{myprocspace}) and (\ref{newspace}) lead to examples of the circle actions 
on the links of affine cones on these singular hypersurfaces which are rationally $S^{2n+1}$ and $S^n\times S^{n+1}$, respectively (cf. Remark \ref{circleaction}).

All known examples of rationally elliptic 
algebraic varieties are rationally connected (which prompted question 1.4 in \cite{amorosbiswas}). The invariants of these singular algebraic varieties 
are in agreement with the Hodge theoretical properties of rationally elliptic singular algebraic varieties discussed in \cite{joinshoji} and \cite{yokura} (cf. Remark  \ref{hodgetate}).

To obtain examples of rationally elliptic varieties with properties mentioned above, we will consider three types of singular spaces, all of which are given by explicit equations:
\begin{equation}\label{kollar}
 H(a_0,\cdots, a_{n+1}): z_0^{a_0}z_1+z_1^{a_1}z_2+\cdots z_n^{a_{n+1}}z_{n+1}+z_{n+1}^{a_{n+1}}z_0=0 \ \ n \ odd
\end{equation} 
\begin{equation}\label{myprocspace}
 V_n^d: \ \  z_1^{d}+z_0z_2^{d-1}+z_2z_3^{d-1}+\cdots z_nz_{n+1}^{d-1}=0
\end{equation} 
\begin{equation}\label{newspace}
W^d_n: z_0z_1^{d-1}+z_1z_2^{d-1}+z_2z_3^{d-1}+\cdots z_nz_{n+1}^{d-1}=0.
\end{equation}
The first class consists of hypersurfaces in a weighted projective space, while the last two are hypersurfaces in $\PP^{n+1}$. 
(\ref{kollar}) and (\ref{myprocspace}) were considered as part of the study of homological projective spaces (cf.\cite{kollaryau}, \cite{proc77} and also \cite{bartheldimca} respectively.).
The third class (\ref{newspace})  is new  \footnote{However, the surface $W_2^3$ coincides with the cubic surface with $D_5$-singularity considered in \cite{brucewall} after change of variables:
$x_0=z_2, x_1=z_1, x_2=z_3,x_3=z_0$. It was also considered in a very different context in  T. Browning and U. Derenthal, Manin's Conjecture for a cubic surface with $D_5$ singularity. Int. Math. Res. Not. IMRN 2009, no. 14, 2620-2647.} and gives singular models of homological even-dimensional quadrics.
We show that these singular spaces have the $\RR$-homotopy type of 
their smooth counterparts, which are known to be rationally elliptic. 

The relevant calculations are presented in the following two sections. 
The topology of the hypersurfaces in the case of families (\ref{myprocspace}) and (\ref{newspace}) 
depend on the explicit form of the characteristic polynomials of the monodromy of the corresponding singularity 
in each of the cases. The calculation was done for (\ref{myprocspace}) in \cite{proc77} and for (\ref{newspace}) is in Section \ref{homquadric}.
In respective sections, we also describe the cohomology algebras as the quotients of polynomial rings (with changed grading in the case 
of $W_n^d$). The cohomology algebras of $V_n^d$ are just the cohomology of twisted projective spaces that appeared in the study of 
complete intersections (cf. \cite{diffstructures}).

Finally, I want to thank S.Yokura for bringing the Hilali conjecture to my attention (cf. \cite{yokura}) which led to collaboration \cite{joinshoji} 
and to my interest in rationally elliptic spaces. I  also want to thank him, J. Amoros as well as the referee of this paper 
for the comments on the various drafts of this note.

\section{Homotopy projective spaces}\label{homproj}

In this section, we show in Theorems \ref{kollarspaces} and \ref{myspaces} that the examples of homology projective spaces from \cite{kollaryau}, \cite{proc77} and \cite{bartheldimca} are rationally elliptic and in final remark point out that some homotopy types of these singular examples can be realized by Fano manifolds. 

\begin{thm}\label{kollarspaces} Let $a_i \in \NN, i=0,\cdots n+1$ where $n \ge 3$ is odd.  Following J.Koll\'ar (cf. \cite{kollaryau}, Sec. 5) 
denote by $d=a_0\cdots a_{n+1}+(-1)^{n+1}$, the determinant of the following system of $n+2$ equations with unknowns $w_i$  
\begin{equation}\label{weightsystem}w_{\bar i}+a_{\overline{i-1}\cdot }w_{\overline{i-1}}=d
\end{equation} 
where for $j\in \ZZ$, $\bar j$ is defined by conditions: $0 \le \bar j \le n+1$ and $j \equiv \bar j \ mod  \ n+2$.  Assume that $$w^*=gcd (w_0,\cdots w_{n+1})=1$$
where $w_i$ are the positive integral solutions to (\ref{weightsystem}). Then the hypersurface $H(a_0, \cdots a_{n+1})$ in 
$\PP(w_0,\cdots, w_{n+1})$ given by the equation (\ref{kollar}):
 $$z_0^{a_0}z_1+z_1^{a_1}z_2+\cdots z_n^{a_{n+1}}z_{n+1}+z_{n+1}^{a_{n+1}}z_0=0$$
 is rationally homotopy equivalent to $\PP^n$. 
In particular, it is rationally elliptic.
\end{thm}

\begin{proof} $H(a_0, \cdots a_{n+1})$ is quasi-smooth (cf. \cite{kollaryau}, Theorem 39 (2)) 
and hence is simply-connected (cf. \cite{dolg} Th. 3.2.4). By \cite{kollaryau} it has  the $\QQ$-homology of 
the projective space. Moreover, the Leray spectral sequence of the fibration with the total space being the link of the singularity 
in $\CC^{n+2}$ given by equation (\ref{kollar}) and the base being the hypersurface $H(a_0,\cdots, a_{n+1})$ shows that 
$H^*(H(a_0,\cdots, a_{n+1}),\QQ)$ is generated by the pullback of  $c_1(\CO_{\PP(w_0,\cdots, w_{n+1})}(1))$. 
The monomorphism of rings $\ZZ[x]/(x^{n+1}) \rightarrow H^*(H(a_0,\cdots, a_{n+1}),\ZZ)$ given by $x \rightarrow   c_1(\CO_{\PP(w_0,\cdots, w_{n+1})}(1))$
has a finite cokernel \footnote{This cokernel is nontrivial, cf. \cite{kawasaki}} since $h^i \ne 0$ for $i \le n$, 
and hence induces an isomorphism of rings over $\QQ$. Moreover, since
$H(a_0,\cdots a_{n+1})$ has only quotient singularities, i.e. is a $\QQ$-manifold, it is formal (cf. \cite{chataur}, Example 3.5) and hence the isomorphism of the cohomology rings induces the isomorphism
of Sullivan's minimal models, which shows the claim. 
\end{proof}

\begin{thm}\label{myspaces} The hypersurface $V_n^d$ ($d,n \in \NN$), given by (\ref{myprocspace}) i.e.  
$$z_1^{d}+z_0z_2^{d-1}+z_2z_3^{d-1}+\cdots z_nz_{n+1}^{d-1}=0$$
has rational homotopy type of $\PP^n$. In particular, it is rationally elliptic. 
It is birational to $\PP^n$. Moreover, the self-intersection of the canonical class $K_{V^d_n}^n$ is positive for $d>n+2$ and can be arbitrarily
large (it is zero for $d=n+1$ and negative otherwise).
\end{thm} 

\begin{proof} The hypersurface (\ref{myprocspace})
is the singular fiber of the family:
\begin{equation}
V^d_n(c): \ z_1^{d}+z_0z_2^{d-1}+z_2z_3^{d-1}+\cdots z_nz_{n+1}^{d-1}=cz_0^d
\end{equation}
has an isolated singularity at $P: z_0=1,z_i=0, i\ge 1$. The Milnor number of this singularity is equal to ${1\over d}[(d-1)^{n+2}+(-1)^n(d-1)]$ and coincides
with the rank of the group of vanishing cycles of a smooth hypersurface in $\PP^n$.  The link of the singularity at $P$ is homeomorphic
to a sphere if $n$ is odd (cf. \cite{proc77}, Theorem 2)). For $n$ even this link is a $(n-1)$-connected $(2n-1)$-manifold $\partial F^d_n$, where $F^d_n$ is the Milnor fiber 
of the singularity at $P$ such that $H_{n-1}(\partial F^d_n,\ZZ)=\ZZ/d\ZZ$ (cf. \cite{proc77} and \cite{top79}).
Moreover, $V_n^d(0) \setminus P$ retracts onto the 
intersection of $V_n^d(0)$ with the hyperplane at infinity $z_0=0$ since the singularity at $P$ is weighted homogeneous. 
One has for $i \ne 0$: 
$$H^i(V^d_n(0),\ZZ)=H^i(V^d_n(0),P,\ZZ)=H_{2n-i}(V^d_n(c),F^d_n,\ZZ).$$
The exact sequence of the 
pair $(V_n^d(c),F_n^d)$ shows that $$H_n(V_0,P,\ZZ)=H_n(V_n^d(c),F^d_n,\ZZ)=H_n(V_n^d(c),\ZZ)/Im(H_n(F^d_n,\ZZ))$$
since $H_{n-1}(F^d_n, \ZZ)=0$. The subgroup $ImH_n(F^d_n,\ZZ) \subset H_n(V_n^d(c),\ZZ)$ is a subgroup of the group of vanishing cycles i.e. 
the orthogonal complement to the class of hyperplane section in $H_n(V_n^d(c),\ZZ)$ and coincides with it since
both subgroups have the same ranks and the same determinants of the intersection form (i.e. $\pm d)$. Clearly 
the quotient of  $H_n(V_n^d(c),\ZZ)$ by the subgroup of vanishing cycles is generated by the coset of the 
class of the $n \over 2$-dimensional linear subspace of smooth hypersurface and hence is $\ZZ$. 
This implies that $V^d_n(0)$ is a $\ZZ$-homology projective space.

Moreover,
\begin{equation}\label{cohomodd} 
H^*(V_n^d(0),\ZZ)=\ZZ[h,y_{n+1}]/h^{n+1},h^{{n+1}\over 2}=dy_{n+1}  \ \ \ n \ \ odd 
\end{equation}
and 
\begin{equation}\label{cohomeven}
H^*(V_n^d(0),\ZZ)=\ZZ[h,y_{n+2}]/h^{n+1},h^{{n\over 2}+1}=dy_{n+2} \ \ \ n  \ even  
\end{equation}
where $y_{n+1} \in H^{n+1}(V^d_n(0),\ZZ)$ for $n$ odd (respectively $y_{n+2} \in H^{n+2}(V^d_n(0),\ZZ)$ for $n$ even)
is the class dual to the linear subspace of dimension ${n-1}\over 2$ for $n$ odd (respectively linear subspace of dimension ${n \over 2}-1$ for $n$ even). \footnote{The spaces with such cohomology ring were called {\it twisted projective spaces} in \cite{diffstructures}}

Indeed, for $n$ odd $V_n^d(0)$ is a topological manifold (since, as was mentioned before, the link of its singularity is homeomorphic to a sphere) 
with generators of $H^{2i}(V_n(0),\ZZ)$ for $i \le n-1$ being $h^i$ where $h$ is the cohomology class 
of the hyperplane section and for $i \ge n+1$ being $y_{n+1}h^{i-n-1}$ where as above $y_{n+1}$ is the class dual to the linear subspace of dimension ${n-1}\over 2$ given in $\PP^{n+1}$
by ${{n+1}\over 2}$ equations: $z_1=0, z_2=z_4=\cdots=z_{n+1}=0$. The multiplication formulas follow from the intersection formulas for the explicit cycles dual to $y_{n+1}$ and $h^i$.

In the case when $n$ is even, the multiplication formulas follow since 
$$H^i(V^d_n(0),P,\ZZ)=H^i(V^d_n(0)\setminus P,\partial F^d_n,\ZZ)=H_{2n-i}(V^d_n(0)\setminus P,\ZZ)=$$
$$H_{2n-i}(V_{n-1}^d(0),\ZZ)=H^{i-2}(V_{n-1}(0),\ZZ)$$
The classes $y_{n}, h^i$ in $H^*(V_{n-1}(0),\ZZ)$ via these isomorphisms define the classes $y_{n+2},h^i$ in $H^*(V_n(0),\ZZ)$ with the relation 
$h^{{n+1} \over 2}=dy_{n+1}$ in $H^n(V_{n-1}(0),\ZZ)$,  becoming the relation 
$h^{{n \over 2}+1}=dy_{n+2}$ in $H^{n+2}(V_n(0),\ZZ)$ where $y_{n+2}$ is the class dual to $\PP^{{n \over 2}-1} \subset V_{n-1}(0) \subset V_n(0)$.

\begin{rem} For $n$ even, the hypersurface $V_n(0)$ contains $n \over 2$-dimensional linear subspace $L_{n\over 2}$ given by $z_1=z_{2i}=0, i=1,\cdots, {n\over 2}$
defining the dual cohomology class in $H^{n\over 2}(V_n(0),{1\over d}\ZZ)$ (the group of coefficients consists of the subgroup of $\QQ$ with 
fractions with denominator being the divisors of $d$). This implies that $H^*(V_n(0),\QQ)$ has presentation in which the relation 
$h^{n\over 2}=dl$ where $l$ is the class dual to $L_{n \over 2}$ replacing the relation in (\ref{cohomeven}) (where $l$ is the class of linear subspace of 
dimension ${n \over 2}-1$). 
\end{rem}

A hypersurface (\ref{myprocspace})
 is simply connected as follows from van Kampen theorem (or from Lefschetz hyperplane section theorem). By \cite{chataur}, Example 4.4, it is formal and its cohomology ring over $\QQ$ is isomorphic to that of 
$\PP^n$ which yields the first claim. The birational map can be obtained using the projection from the singular point $P$ and assigning to a point $Q$ in $z_0=0$ the unique (since the multiplicity of this singular point is $d-1$) intersection point $R(Q)$ of the line 
through this point $Q$ and $P$ such that $R(Q)\ne P$. This correspondence gives a biregular map between Zariski dense subsets of the hyperplane $z_0=0$ and $V_n^d$.
Finally, since the codimension of the singular locus of $V_n^d$ is greater than 2 for $n>3$,  the canonical bundle is given by adjunction formula as $\CO_{V_n^d}(n+2-d)$ and 
the claim about self-intersection follows.
\end{proof}

\begin{rem}\label{fano} The examples of twisted projective spaces that are algebraic varieties in Theorems \ref{kollarspaces} and \ref{myspaces} 
are singular. Smooth examples are given by Fano varieties $X_5, X_{22}$ with Betti numbers $b_2=1,b_3=0$, which are
 twisted projective spaces of degree 5 and 22 respectively. Iskovskikh-Mukai classification shows that these (and $\PP^3,Q_2^3$) are the only 
 rationally elliptic smooth Fano twisted $\PP^3$'s(cf. \cite{morimukai}). \footnote{i.e. $X_5$ should be 
 included in (iii) in Cor. 4.8 in \cite{amorosbiswas}.}
\end{rem}

\section{Homotopy quadrics}\label{homquadric}

In this section, we show that some explicitly given singular hypersurfaces have a $\RR$-homotopy type
of smooth even dimensional quadric. The main step is the calculation of the topological invariants of certain singularities 
in Lemma \ref{keycalc}. 

\begin{thm}\label{homquadricsthm} If $d$ and $n$ are both even, then $W_n^d$ (cf. (\ref{newspace})) given by:
$$ z_1^{d-1}z_0+z_1z_2^{d-1}+z_2z_3^{d-1}+\cdots z_{n}z_{n+1}^{d-1}=0 $$
is a rationally elliptic hypersurface with isolated singularity, $\RR$-homotopy equivalent to a smooth quadric.
The self-intersection of the canonical class $K_{W_n^d}^n$ depending on $d$, can be positive (and arbitrarily large), zero, or negative.
\end{thm}

The main step will be a proof that the cohomology ring of a hypersurface (\ref{newspace}) has the following form:

\begin{dfn} \label{twisted}An $n=2k$-dimensional $a$-twisted quadric is a CW complex $X$ with cohomology ring:
\begin{equation}\label{twistedpresentation} \QQ[h,v]/(hv=0, h^{2k+1}=0, h^{2k}=a v^2),  \ \ \ a\in \QQ \ \ a\ne 0
\end{equation} 
with the following grading: $deg h=2, deg v=2k$.
\end{dfn}

\begin{example}\label{quadric} The cohomology ring of an even dimensional quadric $Q_{2k}$ is given by  (cf. \cite{brieskorn}, Satz 3):
for $k$ even:
\begin{equation} 
H^*(Q_{2k},\ZZ)=\ZZ[h,l]/(h^{k+1}=2hl, l^2=lh^{k},h^{2k+1})
\end{equation}
and for $k$ odd:
$$H^*(Q_{2k},\ZZ)=\ZZ[h,l]/(h^{k+1}=2hl, l^2,h^{2k+1})$$
Here $h\in H^2(Q_{2k},\ZZ)$ is the class of hyperplane section, and $l\in H^{2n}(Q_{2k},\ZZ)$ is the class of linear subspace of dimension $k$.
\footnote{The work \cite{brieskorn} uses generators $e_0,\cdots e_j \cdots e_{2k}, e_k'$ of the groups $e_j \in H^{2j}(Q_{2k},\ZZ)$ which in terms
of $h,l$ are as follows: $e_i=h^i, \ {\rm for}  \  i<k, e_k=l, e_k'=h^k-l, e_i=h^{i-k}l,\ {\rm for} \  i>k$.}
In particular, generators of the integral cohomology groups $H^i(Q_{2k},\ZZ)$ are as follows:
$$\begin{cases} h^i & i < k \\ 
                         h^k, l & i=k \\
                         h^{i-k}l    & k<i \le 2k \\
\end{cases}
$$
The difference in the ring structure in the cases when $k$ is even and odd corresponds to the difference in the self-intersection of subspaces
which is zero or 1 as is seen immediately from the calculation of the Chern class of the normal bundle of $\PP^k \subset Q_{2k}$.
\footnote{Explicitly, the self-intersection of a linear space of dimension $k$ on a smooth $2k$-dimensional hypersurface of degree $d$ 
is given by ${1\over {d}}(1-(1-d)^{k+1})$. This calculation was carried out in (9.4) (ii) \cite{kulkarniwood} (the final expression contains a minor typo).}  

The class of vanishing cycle in the degeneration of smooth quadric $Q_{2k}$ into a cone over a smooth quadric $Q_{2k-1}$ is given by $v=2l-h^k$.
In terms of generators $v,h$ one has 
\begin{equation} H^*(Q_{2k},\QQ)=\QQ[h,v]/(h^{2k+1}=0, hv=0, h^{2k}=(-1)^kv^2).
\end{equation}
In particular, a smooth quadric is a $(-1)^k$-twisted quadric as defined in Definition \ref{twisted}.
\end{example} 

\begin{lemma}\label{keycalc}
(a) The hypersurface $W_n^d$ (as in Theorem \ref{homquadricsthm}) is smooth except for the point $P$ with coordinates $z_i=0, i \ne 0$ where it has an isolated singularity. \footnote{This singularity is not a quotient singularity in general.}

(b) The singularity at $P$ is weighted homogeneous  \footnote{i.e. the equation is a linear combination of monomials $z_1^{i_1}\cdots z_{n+1}^{i_{n+1}}$ such that 
$\sum_{j=1}^{n+1} {i_j \over w_j}=1$.} with weights 
\begin{equation}w_i={{d(d-1)^i} \over {(d-1)^i+(-1)^{i-1}}}, i=1, \cdots n+1
\end{equation}

(c) The characteristic polynomial $\Phi^d_n(t)$ of the local monodromy of singularity of $W^d_n$ satisfies the 
recurrence relations: 
\begin{equation}\label{recurrence}
\Phi^d_{n+1}(t)=\Phi^d_{n-1}(t){{t^{(d-1)^{n+2}}-1}\over {t^{(d-1)^{n+1}}-1} } 
\end{equation}
where $$\Phi^d_0(t)={{t^{d-1}-1}\over {t-1}}$$
In particular, the Milnor number of the singularity at $P$ is 
\begin{equation}\label{milnornumberformula} deg \Phi^d_n(t)={1\over d}[(d-1)^{n+2}-1].
\end{equation}
Moreover, 
\begin{equation}\label{valueat1}
\Phi^d_{2m}(1)=(d-1)^{m+1}
\end{equation}
(d) $W^d_n$ is a rational homology manifold having the $\QQ$-homology of a smooth quadric.
\end{lemma} 

\begin{proof} The local equation of the singularity at $P$ is:
\begin{equation}
z_1^{d-1}+z_1z_2^{d-1}+z_2z_3^{d-1}+\cdots z_nz_{n+1}^{d-1}=0
\end{equation}
Since 
\begin{equation} {{(d-1)^i+(-1)^{i-1}} \over {d(d-1)^i}}+(d-1){{(d-1)^{i+1}+(-1)^{i}} \over {d(d-1)^{i+1}}}=1
\end{equation}
(a) and (b) follow. 

To show (c) we notice that the singularity of $V^d_{n+1}$ can be obtained from $W_n^d$ by adding $d$-th power of 
an independent variable. Denoting, by $\Delta_{n+1}(t)$ (resp. $\Phi_n(t)$) the characteristic polynomial of the monodromy of singularity of 
$V_{n+1}^d$ (resp. $W_n^d$) it follows from a theorem of Thom and Sebastiani (cf. \cite{thomseb}) that
\begin{equation}\label{charpolrel1}
\Delta_{n+1}(t)=\prod_{i=1}^{d-1} \Phi^d_n(\omega^i_dt)
\end{equation} 
($\omega_d$ is a primitive root of unity of degree $d$). This implies that for $j=0,\cdots, d-1$ one has 
\begin{equation}\label{carpolrel2}
\Delta_{n+1}(\omega_d^jt)=\prod_{i=1}^{d-1} \Phi^d_n(\omega^{i+j}_dt)=
\prod_{i=0, i\not\equiv -j (mod \ d)}^{d-1} \Phi^d_n(\omega^i_dt)
\end{equation}
Multiplying the relations (\ref{carpolrel2}) for $j=0, \cdots d-1$ we obtain: 
$$ \prod_{j=0}^{d-1} \Delta_{n+1}(\omega_d^jt)=\prod_{j=0}^d \prod^{d-1}_{i=0, i\not\equiv -j({\rm mod} \ d)} \Phi^d_n(\omega_d^it)=(\Phi^d_n(t)\prod_{i=1}^{d-1} \Phi^d_n(\omega_d^it))^{d-1}.
$$
Therefore
\begin{equation}\label{carpolrel}\Phi^d_n(t)={{[\prod_{i=0}^{d-1} \Delta_{n+1}(\omega_d^it)]^{1\over {d-1}}} \over {\Delta_{n+1}(t)}}
\end{equation}
Recall from \cite{proc77} Theorem 2 (iii) that 
\begin{equation}
\Delta_{n+2}(t)=\Delta_n(t) {{(t^{d(d-1)^{n+2}}-1)(t^{(d-1)^{n+1}}-1)}\over {(t^{(d-1)^{n+2}}-1)(t^{d(d-1)^{n+1}}-1)} },  \ \ \ \ \Delta_1(t)={{(t^{d(d-1)}-1)(t-1)}\over {(t^{(d-1)}-1)(t^{d}-1)} } 
\end{equation}
This and (\ref{carpolrel}) imply that 
$$\Phi^d_{n+1}=\left( \prod_{i=0}^{d-1} \Delta_n(\omega_d^it)\right)^{1\over {d-1}} \left[\prod_{i=0}^{d-1} {{((\omega_d^it)^{d(d-1)^{n+2}}-1)((\omega_d^it)^{(d-1)^{n+1}}-1)}\over {((\omega_d^it)^{(d-1)^{n+2}}-1)((\omega_d^it)^{d(d-1)^{n+1}}-1)} }\right]^{1\over {d-1}} \times $$
$${{(t^{(d-1)^{n+2}}-1)(t^{d(d-1)^{n+1}}-1)}\over {\Delta_{n}(t)(t^{d(d-1)^{n+2}}-1)(t^{(d-1)^{n+1}}-1)}}=$$
$$\Phi^d_{n-1}(t)\left [\left({{t^{d(d-1)^{n+2}}-1}\over {t^{d(d-1)^{n+1}}-1} }\right)^d{{(t^{d(d-1)^{n+1}}-1)}\over {(t^{d(d-1)^{n+2}}-1)} }\right]^{1\over {d-1}}
{{t^{(d-1)^{n+2}}-1}\over {t^{d(d-1)^{n+2}}-1}} \cdot
{{t^{d(d-1)^{n+1}}-1}\over {t^{(d-1)^{n+1}}-1}}
=
$$
$$\Phi^d_{ n-1}(t){{t^{(d-1)^{n+2}}-1}\over {t^{(d-1)^{n+1}}-1} }
$$
and similarly the expression for $\Phi^d_0(t)$ obtained from recurrence (\ref{carpolrel}). 
The recurrence relation (\ref{recurrence}) shows that $deg \Phi_{n+1}^d=deg \Phi_{n-1}^d+(d-2)(d-1)^{n+1}$ and (\ref{milnornumberformula})
follows, while (\ref{valueat1}) is clear.
This shows (c). 
Note that the calculation of $\Phi^d_n(t)$ can be obtained directly from results of \cite{milnororlik} similarly
to the calculation in \cite{proc77}.

To show (d), denote by $W_n^d(c)$ the smoothing
\begin{equation}\label{singularityP}
  z_0z_1^{d-1}+z_1z_2^{d-1}+z_2z_3^{d-1}+\cdots z_nz_{n+1}^{d-1}=cz_0^d
\end{equation}
of $W_n^d$ ($\vert c \vert$ is sufficiently small and let $M_n^d\subset W_n^d(c)$ be the Milnor fiber of 
singularity (\ref{singularityP}).
The exact sequence of pair $(W_n^d(c),M_n^d)$, since $W_n^d(0)=W_n^d(c)/M_n^d$,  shows that 
\begin{equation}\label{homologyquadric}
H^i(W_n^d,\ZZ)=H^i(W_n^d(c),\ZZ)=
\begin{cases} 0 & i \ is \ odd \\ 
   \ZZ & i \ is \ even \  i \ne n  
\end{cases}
\end{equation} 
Moreover, 
\begin{equation}
 rk H^n(W_n^d)=rk H^n(W_n^d(c))-rk H^n(M_n^d)
\end{equation}
The subspace of vanishing cycles $Van(H^n(M_n^d(c),\QQ)$ has codimension 1 in $H_n(W_n^d(c),\QQ)$ and 
contains $H_n(M_n^d,\QQ)$ as a proper subset. One has 
 $$rk Van(H^n(M_n^d(c),\QQ)={1\over d}[(d-1)^{n+2}+(d-1)]$$  
 Hence (\ref{milnornumberformula}) shows that 
 $H_n(M_n^d,\QQ)$ has codimension 1 in $Van(H^n(M_n^d(c),\QQ)$ and therefore 
 $rk H^n(W_n^d,\QQ)$ coincides with the rank of middle-dimensional homology of an even-dimensional quadric  (cf. Example 
 \ref{quadric}). Since $\Phi_{2m}^d (t)\ne 0$ (cf. \ref{valueat1}), it follows that the link of the only singularity of $W_n^d$ is a $\QQ$-homology 
 sphere which shows d). 
\end{proof} 
\begin{prop} Two CW-complexes that are
  formal and are $a$-twisted even-dimensional quadrics have the same real homotopy type if and only if the signs of corresponding rationals $a$ are
  the same. 
\end{prop} 
\begin{proof} Indeed, denoting by $R_a$ the graded algebra given by (\ref{twistedpresentation}) we see that the 
map $h \rightarrow h, v \rightarrow {\vert a \vert}^{1 \over 2}v$ gives an isomorphism between $R_{sign(a)}$ 
and $R_a$. Clearly $R_1$ and $R_{-1}$ are not isomorphic over $\RR$ 
(though they are isomorphic over $\CC$). Since for formal spaces the cohomology algebra determines the minimal 
model, the claim follows.
\end{proof}
\begin{rem} The same argument shows that such CW-complexes having homological dimension $n=2k$ are rationally homotopy equivalent if and only if 
the quadratic forms on middle-dimensional cohomology given by $(h^k,v)=h^{2k}-av^2$ are $\QQ$ equivalent which
is described by the class of $a$ in $\QQ^*/(\QQ^*)^2$. In particular, the twisted quadrics $W_{2k}^d$ are not rationally 
homotopy equivalent. This is reflection of the phenomenon discussed in \cite{sczarba} (\cite{herrmann} has explicit examples)
\end{rem}

\begin{prop}\label{typeofquadric} $n=2k$-dimensional hypersurface $W^d_{2k}$ is $\RR$-homotopy equivalent to the $(-1)^{k}$-twisted quadric.
\end{prop}

\begin{proof} It follows from (\ref{homologyquadric}) that the cohomology ring $H^*(W_{2k}^d,\RR)$ is generated by the
degree 2 class of hyperplane section $h$ and by a class $v$ of degree $2k$ that is image of a vanishing cycle of degeneration 
of smooth hypersurface $W_{2k}^d(c)$ in its degeneration into $V_{2k}^d$. One has $h^i \ne 0, i\le 2n$ and $hv=0$ 
and $h^{2n}=av^2, a \ne 0$ since $H^{2k}(W_{2k}^d)=\RR$.
 
J.Steenbrink showed (cf. \cite{steenbrink} Prop. 4.14) that $\mu-\mu_{(-1)^k}$ is even, where $\mu$ is the Milnor number and $\mu_{(-1)^k}$ is the (positive or negative) inertia index of of the intersection form
on the middle cohomology of Milnor fiber of a function in $2k+1$ variables. 
In particular, the parity of $\mu-\mu_{(-1)^k}$ is given by the dimension modulo 4 and is independent of $d$. 
Since $\Phi^d_{2k}(1)\ne 0 $ and $\Delta_{2k}(1)\ne 0$ 
it follows that the inertia index $\mu_0$ is zero for singularities of both hypersurfaces $W_{2k}^d$ and $V_{2k}^d$ and hence $\mu-\mu_{(-1)^k}=\mu_{(-1)^{k-1}}$ is even in both cases. Since ${\rm rk} H^k(V_{2k}^d,\RR)/H^n(W_{2k}^d,\RR)=1$ it follows that  
the restriction of the intersection form to $v \in H^k(V_{2k}^d,\RR)/H^n(W_{2k}^d,\RR)$ has sign $(-1)^k$ and hence the signature 
of the intersection form on $H^n(W_{2k}^d,\RR)$ is the same as the signature of a smooth quadric. Therefore the cohomology 
graded algebras of $H^*(W_{2k}^d,\RR)$ and smooth quadric are isomorphic. Formality (cf. \cite{chataur}) implies that this isomorphism is induced
by a real homotopy equivalence. This also implies all the claims on the homotopy type in Theorem \ref{homquadricsthm}.
\end{proof}

\begin{prop} Odd dimensional quadircs  $Q_{2k+1}$ are rationally homotopy equivalent to $\PP^{2k+1}$.  
\end{prop} 

\begin{proof} Since a quadric  $Q_{2k+1}$ is a ramified double cover of $\PP^{2k+1}$, the induced 
ring homomorphism is an isomorphism over $\QQ$ since the ranks of cohomology groups are the same.
\end{proof}

\begin{rem}\label{circleaction} The affine cones over (\ref{myprocspace}) and (\ref{newspace}) (i.e. hypersurfaces in 
$\CC^{n+2}$ given by the same equations as in projective space) have 1-dimensional singular locus and have as their links $L(V_n^d)$ 
and $L(W_n^d)$ respectively, that are  $2n+1$ dimensional  $\QQ$-manifolds that are $\QQ$-homology spheres in the case (\ref{myprocspace}) and $H_n(L(W_n^d,\QQ))=H_{n+1}(L(W_n^d,\QQ))
=\QQ$ (i.e. a rational $S^n\times S^{n+1}$).
\end{rem}

\begin{rem}\label{hodgetate} The invariants of singular algebraic varieties that are rationally elliptic and additional structures on them were 
considered in \cite{joinshoji} and \cite{yokura}. 
The examples in this note all have the mixed Hodge structure that is Hodge-Tate 
 (the only type afforded by the structure of their cohomology groups) in agreement with (\cite{yokura}, Conj.10.4) and other 
examples of quasi-projective varieties (including singular, e.g. \cite{wiemeler}). 
\end{rem}

\begin{rem}\label{affine} Let $H_{1}$ be the hyperplane $z_{1}=0$. Then $W^d_n\setminus (W_n^d\cap H_1)$ is 
a smooth affine hypersurface bihomolomorphic to $\CC^n$. The identification is given by the 
map: $(z_2,\cdots,z_{n+1}) \rightarrow (-(z_2^{d-1}+z_2z_3^{d-1}+\cdots z_nz_{n+1}^{d-1}), z_2,\cdots, z_{n+1})$
(which is the inverse of the restriction of projection $(z_0,z_2,\cdots, z_{n+1})$ to $W^d_n\setminus (W^d_n \cap H_1)$).
\end{rem}

\section{Appendix: Smooth rationally elliptic algebraic varieties in dimensions 2 and 3.}

We show that the rationally elliptic Fano 3-folds listed in \cite{amorosbiswas} (and Fano manifold $V_5$)
give a complete list of rationally elliptic smooth threefolds that are birationally minimal (i.e. the canonical class is nef or are Fano threefolds). 
This implies the finiteness of deformation types of such rationally elliptic threefolds. On the other hand, the set of projective line bundles 
over $\PP^2$ form an infinite family of non-diffeomorphic rationally elliptic 3-folds having the same real homotopy type.

\begin{thm}\label{lowdim}
1. The types of rationally elliptic algebraic surfaces up to biregular equivalent are as follows:
1) $\PP^2$, or 2) Hirzebruch surfaces $H_n, n \ge 0$ or 3) a fake quadric i.e. 
a surface of general type with $p_g=q=0, K^2=8, c_2=4$. Two Hirzebruch surfaces $H_n, H_m$ are diffeomorphic if and only if $m=n \ mod \ 2$.

2. a) There are no rationally elliptic K\"ahler three-folds with nef canonical bundle.
The  3-folds with nef anti-canonical bundles that are rationally elliptic have a finite number of 
deformation types (and hence finitely many integral homotopy types or diffeomorphism types).

b) The 3-folds $\PP^2(E)$, where $E$ is a holomorphic bundle of rank 2, form an infinite family of manifolds having the same real homotopy type but 
pairwise homotopically distinct.  More precisely, any such a 3-fold is isomorphic to $F_n=\PP^2({\cal O}\oplus {\cal O}(n)), n 
\ge 0$ and 
$F_n$ is homotopy equivalent to $F_m$ if and only if $m=n$.  
\end{thm}

\begin{proof}
The case of surfaces is a restatement of Theorem 1.1 in \cite{amorosbiswas} and well-known facts about Hirzebruch surfaces 
 (cf. \cite{fm}. Ch.I sect.1.2).

If $X$ is a 3-fold with nef canonical bundle, then $c_1^3 \le 0$. It follows from Myiaoka's inequality (cf. \cite{miya}) that $c_1c_2\le 0$.
But the Riemann-Roch and the structure of the Hodge diamonds in \cite{amorosbiswas} (see also \cite{herrmann}), 
since $h^{1,0}=h^{2,0}=h^{3,0}=0$, show that 
 $c_1c_2=24\chi({\cal O}_X)$ is 24. Hence $X$ 
is rationally elliptic if and only if $X$ admits a contraction of anticanonical divisor which is nef. In the latter case i.e.  the Fano case, smooth rationally elliptic 
 varieties were identified in \cite{amorosbiswas} (and $X_5$, cf. Remark \ref{fano}) and hence belong to finitely many deformation types.  

To describe the topology of the Mori fiber space $\PP(E)$, recall that the classification of Wall (\cite{wall}) (completed in \cite{jupp}, \cite{zubr}) 
shows that the diffeomorphism type of a $C^{\infty}$-manifold $X$ having real dimension 
6, is determined  by the Betti number $b_3(X)$, the trilinear form on $H^2(X,\ZZ)$ given by the cup product, the second Stiefel-Whitney class 
$w_2 \in H^2(X,\ZZ/2\ZZ)$ and
 a class in $H^4(X,\ZZ)=Hom(H^2(X,\ZZ),\ZZ)$ 
representing the Pontryagin class $p_1=c_1^2-2c_2 $. We can assume that $E=\CO_{\PP^2}\oplus \CO_{\PP^2}(n), n\in \ZZ^{\ge 0}$. It follows 
that $H^*(\PP^2(E),\ZZ)=\ZZ[x,y]/(x^3, y^2+nxy)$. 
The oriented homotopy type of a simply-connected real 6-fold is given by $H^2(X,\ZZ),\mu(x)=x\cap x \cap x[X] \in \ZZ, x\in H^2(X,\ZZ), p_1 \ mod \ 48$
(the latter is an invariant of the homotopy type). In the case $X=\PP(E)$ and the basis $x,y$ above,  
$\mu(ax+by)=3a^2b-3nab^2+n^2b^3, a,b \in \ZZ$. Equivalent binary cubic forms have 
equivalent discriminants (the latter for $aX^3+bX^2Y+cXY^2+dY^3$ is given by 
$\Delta=b^2c^2+18abcd-4ac^3-4b^3d-27a^2d^2$ cf. \cite{hambleton}) and since $\PP^2(\CO\oplus \CO(n))=\PP^2(\CO(-n)\otimes (\CO\oplus \CO(n))$
the claim follows as the direct calculation shows that the discriminant is determined by $n^2$.
\end{proof}

\end{document}